\documentclass[12pt]{article}

\usepackage{amssymb}
\usepackage{amsmath}
\usepackage{euscript}

\begin{document}

\medskip
\centerline{\large\bf Final Dynamics of Systems of Nonlinear}

\centerline{\large\bf Parabolic Equations on the Circle}

\bigskip
\centerline{\large\bf A. V. Romanov}

\bigskip
\noindent \textit{Abstract.} We consider the class of dissipative
reaction-diffusion-convection systems on the circle and obtain
conditions under which the final (at large times) phase dynamics
of a system can be described by an ODE with Lipschitz vector field
in $\mathbb{R}^{N}$. Precisely in this class, the first example of
a parabolic problem of mathematical physics without the indicated
property was recently constructed.

\bigskip



\medskip
\centerline{\large\bf 1. Introduction}
\medskip

The problem of describing the final (at large times) dynamics of
dissipative semilinear parabolic equations (SPE)
$$
\partial_{t}u=G(u)
\eqno (*)
$$
(see~[5]) with a Hilbert phase space~$\mathcal{X}$ by ordinary
differential equations (ODE) in~${\mathbb R}^{N}$ has been
attracting researcher's attraction for a long time. In fact, it is
required to separate finitely many ``determining'' degrees of
freedom of an infinite-dimensional dynamical system. In this case,
the key geometric object is the so-called (global)
\textit{attractor} [1,13,16], i.e., the connected compact
invariant set $\mathcal{A}\subset \mathcal{X}$ that uniformly
attracts bounded subsets $\mathcal{X}$ as $t\to +\infty $.

The required ODE can sometimes be implemented as an
\textit{inertial form} [13,16,19] obtained by restricting the
initial equation to an \textit{inertial manifold}, i.e, a
finite-dimensional invariant $C^{1}$-surface ${\mathcal{M}}\subset
\mathcal{X}$ containing the attractor and exponentially attracting
(with asymptotic phase) all trajectories of~$(*)$ as $t\rightarrow
+\infty$. The theory of inertial manifolds originally encountered
systematic difficulties, and several alternative concepts of
\textit{finite-dimensional reduction} of SPE have therefore been
developed starting from [3,12,14,15]. Following~[14], we will say
that the dynamics of $(*)$ on the attractor (final dynamics) is
finite-dimensional if there exists an ODE in~${\mathbb R}^{N}$
with Lipschitz vector field, resolving flow $\{ \, \Theta_{t}
\}_{t\in \mathbb{R}}$, and invariant compact set ${\EuScript
K}\subset \mathbb{\mathbb{R}^{N}}$ such that the phase semiflows
$\{\, \Phi_{t}\}_{t\geq 0}$ of equation~$(*)$ on ${\EuScript A}\;$
and $\{ \, \Theta_{t} \}_{t\geq 0}$ are Lipschitz conjugate on
${\EuScript K}$. The existence of the inertial manifold implies
that the dynamics is finite-dimensional on the attractor and, in
general, looks like a more attractive property. Indeed, in the
first case, the inertial form provides an exponential asymptotics
of any solution of the equation at large times, and in the second
case, we have an ODE reproducing the original dynamics only on the
attractor itself. Nevertheless, the fact that the dynamics is
finite-dimensional on ${\EuScript A}$ means that the structure of
limit regimes of SPE with infinitely many degrees of freedom is no
more complicated than the structure of similar regimes of an ODE
with Lipschitz vector field in~${\mathbb R}^{N}$.

In this paper, we consider the problem of whether the final
dynamics is finite-dimensional for 1D systems of
reaction-diffusion-convection equations
$$
\partial_{t}u=D\partial_{xx}u-u +f(x,u)\partial_{x}u+g(x,u),
\eqno (1.1)
$$
where $u=(u_{1},\dots,u_{m})$ and $f$ and $g$ are sufficiently
smooth matrix and vector functions.We assume that $x\in J$, where
$J$ is a circle of length~1. The matrix of diffusion
coefficients~$D$ is assumed to be diagonal,
$D=\mathrm{diag}\{d_{j}\},\,d_{j}>0$. As the phase space we choose
an appropriate space $\mathcal{X}\subset C^{1}(J,\mathbb{R}^{m})$
in the Hilbert semiscale $\{X^{\alpha}\}_{\alpha\geq 0}$ generated
by a linear positive definite operator $u\rightarrow u-Du_{xx}$ in
$X=L^{2}(J,\mathbb{R}^{m})$. We postulate that evolution equation
(1.1) is dissipative in~$\mathcal{X}$ and there exists the
attractor $\mathcal{A}\subset \mathcal{X}$ consisting of functions
$u=u(x),\; u\in C^{1}(J,\mathbb{R}^{m})$. The algebraic structure
of the ``convection matrix'' $f=f(x,u),\;f=\{f_{ij}\}$,
$i,j\in\overline{1,m}$, on the convex hull
$\mathrm{co}\,\mathcal{A}\subset \mathcal{X}$ plays an important
role. We will highlight the case of the scalar diffusion matrix
$D=dE $, where $d =\mathrm{const}$ and $E$ is the identity matrix.

For scalar equations of the form (1.1), the fact that the dynamics
is finite-dimensional on the attractor was established in~[15]. In
the vector case, the final dynamics of systems (1.1) with scalar
diffusion matrix~$D$ and spatially homogeneous nonlinearity
$f(u)\partial_{x}u+g(u)$ was studied in~[8], and the second
restriction seems to be technical. The existence of an inertial
manifold was proved in~[8] for the scalar equation ($m=1$), and
for $m>1$, it was proved under the assumption that the function
matrix $f(u)$ is diagonal with a unique nonzero element in a
convex neighborhood of the attractor. The results obtained in~[8]
are based on a non-local change of the phase variable~$u$ which
``decreases'' the dependence of the nonlinear part (1.1)
on~$\partial_{x}u$ and allows using the well-known ``spectral gap
condition``.

Generalizing and developing the approach in~[15], we study whether
the dynamics is finite-dimensional on the attractor, but we do not
consider the problem of existence of an inertial manifold for
systems of periodic equations (1.1). At the same time, we here
consider the case of \textit{nonscalar} diffusion matrix~$D$. We
prove that the limit dynamics is finite-dimensional for wide
classes of systems (1.1). Now, omitting the details related to the
choice of phase space and dissipativity conditions, we formulate
the main results of the paper as follows.

 \textit{The phase dynamics on the attractor of system $(1.1)$
is finite-dimensional if any of the following three conditions is
satisfied}.

(A) \textit{The convection matrix $f=\mathrm{diag}$ on ${\rm
co}\,\mathcal{A}$} (Theorem~4.3).

(B) \textit{The diffusion matrix $D$ is scalar.  For all $(x,u)\in
J\times \mathrm{co}\,{\mathcal{A}}$, the numerical matrices
$f(x,u(x))$ have $m$ distinct real eigenvalues and commute with
each other} (Theorem~4.5).

(C) \textit{The diffusion matrix $D$ is scalar.  For all $(x,u)\in
J\times \mathrm{co}\,{\mathcal{A}}$, the matrices $f(x,u)$ are
symmetric and commute with each other} (Theorem~4.6).

In the case~(A), we have $Df=fD$ on $\mathrm{co}\,\mathcal{A}$.
The assumptions that the matrices are commutative can
conditionally be formulated as the \textit{consistency of
convection with diffusion} and the \textit{self-consistency of
convection} on the convex hull of the attractor. Usually, the
attractor $\mathcal{A}$ of system (1.1) can be localized in a ball
$B\subset \mathcal{X}$ centered at zero. Since the embedding
$\mathcal{X}\rightarrow C(J,\mathbb{R}^{m})$ is continuous, it is
actually sufficient to verify the conditions on~$f=f(x,u)$ in
assertions~(A), (B), and~(C) for $x\in J$, $u\in
\mathbb{R}^{m}:|u|<r$ with an appropriate $r>0$.

In the class of one-dimensional systems (1.1), was constructed [8,
Theorem~1.2] the first example of semilinear parabolic equation of
mathematical physics (actually, a system of eight equations with
scalar diffusion) that does not demonstrate any finite-dimensional
dynamics on the attractor. This class seems to be a good testing
ground for understanding where the finite-dimensional final
dynamics of semilinear parabolic equations terminates and the
infinite-dimensional final dynamics begins.

The results of the paper can be generalized to systems on the
circle of the form
$$
\partial_{t} u=D\partial_{xx}u+f(x,u,\partial_{x}u)
\eqno (1.2)
$$
with a smooth vector function $f=(f_{1},\dots,f_{m})$. Such
systems with various boundary conditions can be reduced (see~[7,
8]) to the form (1.1) by the termwise differentiation and an
appropriate change of the variable. The fact that the final
dynamics is finite-dimensional for scalar equations (1.2) was
already proved in~[15].

We here do not consider the Dirichlet and Neumann boundary
conditions for systems of the form~(1.1) on $(0,1)$, this can be
studied in a subsequent paper. The existence of an inertial
manifold is proved in a similar situation in~[7] for systems of
general form~(1.2) with $f=f(u,u_{x})$ and a scalar diffusion
matrix.

The paper is organized as follows. Section~2 contains necessary
information about abstract SPE and the conditions for their final
dynamics to be finite-dimensional. In Section~3, it is shown how
these conditions can be applied to parabolic systems~(1.1). The
main results are obtained in Section~4. In the short Section~5, we
present several examples of system (1.1) which admit a
finite-dimensional final dynamics. Finally, in Section~6, we
discuss alternative approaches to the problem of
finite-dimensional reduction of systems~(1.1).

\bigskip
\centerline{\large\bf 2. General information}
\medskip

First, we consider the abstract dissipative SPE
$$
\partial _{t} u=-Au+F(u)
\eqno (2.1)
$$
in a \textit{real} separable Hilbert space~$X$ with scalar product
$(\cdot,\cdot)$ and the norm $\|\cdot\|$. We assume that the the
unbounded positive definite linear operator~$A$ with domain of
definition $\mathcal{D}(A)\subset X$ has a compact resolvent. We
assume that $X^{\alpha} =\mathcal{D}(A^{\alpha})$ with $\alpha \ge
0$. Then $\left\|u\right\|_{\alpha}=\left\|A^{\alpha}u\right\|$,
$X^{0}=X$, and $X^{1}=\mathcal{D}(A)$. For arbitrary Banach spaces
$Y_{1}$ and $Y_{2}$, we let $BC^{\nu}(Y_{1},Y_{2})$,
$\nu\in\mathbb{N}_{0}$, denote the class of $C^{\nu}$-smooth
mappings $Y_{1}\rightarrow Y_{2}$ that are bounded on balls. We
assume that a nonlinear function $F$ belongs to
$BC^{2}(X^{\alpha},X)$ for some $\alpha \in[0,1)$ and equation
(2.1) is \textit{dissiparive}, i.e., generates a resolving
semiflow $\{\Phi_{t}\}_{t\ge 0}$ in the phase space $X^{\alpha }$
and there exists a \textit{retracting ball}
$\mathcal{B}_{a}=\{u\in X^{\alpha}:\|u\|_{\alpha}<a\}$ such that
$\Phi_{t}\mathcal{B}_{r}\subset \mathcal{B}_{a}$ for any ball
$\mathcal{B}_{r}:\|u\|_{\alpha}<r$ for $t>t^{*}(r)$. In this case,
the semiflow $\{\Phi_{t}\}$ inherits~[5] the $C^{2}$-smoothness,
and there exists the \textit{compact attractor}
$\mathcal{A}\subset\mathcal{B}_{a}$ consisting of all bounded
complete trajectories $\{u(t)\}_{t\in \mathbb{R}}\subset X^{\alpha
}$ and uniformly attracting balls $X^{\alpha}$ as $t\rightarrow
+\infty$. In fact, $\mathcal{A}\subset X^{1}$ due to the
\textit{smoothing action} of the parabolic equation [5,
Section~3.5].

The embeddings $X^{\sigma}\subset X^{\alpha}$ with $\alpha<\sigma<1$
are dense and compact, and ${\|u\|_{\alpha}\leq c\|u\|_{\sigma}}\,$,
$c=c(\alpha,\sigma)$, for $u\in X^{\sigma}$.
Moreover, the proof of Theorem~3.3.6 in~[5] can be used to derive the estimate
$\|\Phi_{1}u\|_{\sigma}\le L(r)\|u\|_{\alpha}$ on the balls
$\mathcal{B}_{r}\subset X^{\alpha}$.
This implies that $F\in BC^{\nu}(X^{\sigma},X)$ if $F\in BC^{\nu}(X^{\alpha},X)$
and the $X^{\alpha}$-dissipativity implies the $X^{\sigma}$-dissipativity.
Thus, in all constructions related to SPE (2.1), one can replace
the nonlinearity index~$\alpha$ with any value $\sigma\in(\alpha,1)$.
The linear operator $A:X^{\vartheta+1}\rightarrow X^{\vartheta}$
is positive definite in $X^{\vartheta}$ with $\vartheta>0$.
If $F\in BC^{2}(X^{\vartheta+\alpha},X^{\vartheta})$,
then one can consider (2.1) in the pair of spaces
$(X^{\vartheta},X^{\vartheta+\alpha})$ instead of $(X,X^{\alpha})$.
In this case, the phase dynamics preserves all its properties listed above.

We say that the phase dynamics of (2.1) is asymptotically
finite-dimensional if there exists an \textit{inertial manifold},
i.e., a smooth finite-dimensional invariant surface
$\mathcal{M}\subset X^{\alpha}$ containing the attractor and
exponentially attracting (with asymptotic phase) all solutions
$u(t)$ at large times. Such a manifold is usually [13,16,19] a
Lipschitz graph over the highest modes of the operator~$A$. The
restriction of SPE (2.1) to $\mathcal{M}$ is an ODE in ${\mathbb
R}^{N},\;N=\dim\mathcal{M}$ which completely describes the final
dynamics of the original evolution system.

A less rigorous approach to the problem of finite-dimensional
limit dynamic of SPE was proposed in~[14,15]. So the dynamics of
(2.1) on the attractor is finite-dimensional if, for some ODE
$\partial_{t}x=h(x)$ in ${\mathbb R}^{N}$ with $h\in {\rm
Lip}({\mathbb R}^{N},{\mathbb R}^{N})$ and resolving flow
$\{\Theta_{t}\}_{t\in {\mathbb R}}$, there exists an invariant
compact set $\mathcal{K}\subset{\mathbb R}^{N}$ such that the
dynamical systems $\{\Phi_{t}\}$ on $\mathcal{A}$ and
$\{\Theta_{t}\}$ on $\mathcal{K}$ are Lipschitz conjugate for
$t\ge 0$. The properties of the dynamics to be asymptotically
finite-dimensional and to be finite-dimensional on the attractor
have not yet been separated; there is a hypothesis~[19] that they
are equivalent.

Here are two criteria for the dynamics to be finite-dimensional on
the attractor~[14] under the assumption that $F\in
BC^{2}(X^{\alpha},X)$.

(Fl) The phase semiflow on $\mathcal{A}$ can be extended to the Lipschitz flow:
$$
\left\|\Phi_{t}(u)-\Phi_{t}(v)\right\|_{\alpha}
\le M\left\|u-v\right\|_{\alpha} e^{\kappa|t|},\quad t\in {\mathbb R},
$$
where $M>0$ and $\kappa\geq 0$ depend only on~$\mathcal{A}$.

(GrF) The attractor is a Lipschitz graph over the lowest Fourier
modes:
$$
\left\| Pu-Pv \right\|_{\alpha}\ge M \left\|u-v\right\|_{\alpha},\quad M=M(\mathcal{A}),
$$
for some finite-dimensional spectral projection $P\in\mathcal{L}(X^{\alpha })$
of the operator~$A$ and all $u,v\in\mathcal{A}$.

Property (GrF) was established for scalar equation (1.1) in~[9]
independently of the results obtained in~[14,15]. We shall further
use other sufficient conditions for the dynamics to be
finite-dimensional on the attractor, which were obtained
in~[15]\footnote{By misunderstanding, an important assumption
that~$X$ is real was not mentioned in~[15].}. Assume that
$G(u)=F(u)-Au\,$ is the vector field of (2.1),
$\mathcal{N}=\mathcal{A}\times\mathcal{A}\subset X^{\alpha}\times
X^{\alpha}$ is a compact set, and $Y$ is a Banach space.
\medskip

\textbf{Definition 2.1 ([15]).} A continuous field
$\Pi:\mathcal{N}\rightarrow Y$ is said to be regular if, for any
$u,v\in \mathcal{A}$, the function
$\Pi(\Phi_{t}u,\Phi_{t}v):[0,+\infty)\rightarrow Y$ belongs to the
class~$C^{1}$ and its derivative $\partial_{t}\Pi(u,v)$ at zero is
bounded uniformly with respect to $(u,v)\in \mathcal{N}$.

The smoothness of the semiflow $\{\Phi_{t}\}$ and the invariancy
of the compact set $\mathcal{A}\subset X^{\alpha}$ imply the
regularity of the identical embedding $\mathcal{N}\rightarrow
X^{\alpha}\times X^{\alpha}$ and hence the regularity of any field
$\Pi:\mathcal{N}\rightarrow Y$ that can be continued to a
$C^{1}$-mapping into the $(X^{\alpha}\times
X^{\alpha})$-neighborhood of the set $\mathcal{N}$. In this
situation, $\partial_{t}\Pi(u,v)=\mathrm{D}\Pi(u,v)(G(u),G(v))$,
where $\mathrm{D}$ is the Frechet differentiation. The regular
fields $\Pi:\mathcal{N}\rightarrow Y$ form a linear structure
which is also multiplicative if $Y$ is a Banach algebra. In the
last case, if the elements of $\Pi(u,v)\in Y$ are invertible, then
the field $\Pi^{-1}$ is also regular, and
$\partial_{t}\Pi^{-1}=-\Pi^{-1}(\partial_{t}\Pi)\Pi^{-1}$ for
$(u,v)\in \mathcal{N}$. We start from the decomposition
$$
G(u)-G(v)=(T_{0}(u,v)-T(u,v))(u-v), \quad (u,v) \in \mathcal{\mathcal{N}},
\eqno (2.2)
$$
of the vector field $G(u)$ on $\mathcal{A}$, where
$T_{0}\in\mathcal{L}(X^{\alpha})$ and $T\in \mathcal{L}(X^{1},X)$
are unbounded linear sectorial operators in~$X$ similar to normal
ones. We write
$$
\Gamma _{a}=\{z\in {\mathbb C}:\mathrm{Re}\, z=a\},\quad
\Gamma (a,\xi )= \{ z\in {\mathbb C}:a-\xi \le \mathrm{Re}\, z\le a+\xi \}
$$
for $a>\xi>0$ and assume that, for some $c>0$, $\theta \in [0,1]$,
the \textit{total spectrum}
$$
\Sigma_{T} =\bigcup_{u,v\in \mathcal{A}} \mathrm{spec}\,T(u,v)
$$
is localized in the domain
$$
\Omega (c,\theta )=\{ x+iy\in{\mathbb C}:\left|y\right|<cx^{\theta }\},\quad x>0\,.
\eqno (2.3)
$$
Let $\beta =\alpha /2$ for $0\le \theta \le \alpha /2$, and let
$\beta =(\alpha +\theta )/3$ for $\alpha /2<\theta \le 1$. Assume
that the set $\mathbb C\backslash \Sigma_{T}$ contains strips
$\Gamma (a_{k},\xi_{k})$ with $a_{k},\xi_{k} \to \infty$ as $k \to
+\infty $.
\medskip

\textbf{Theorem 2.2} (see [15, Theorem~2.8]). \textit{Assume that}
$$
T(u,v)= S^{-1}(u,v)H(u,v)S(u,v)
\eqno (2.4)
$$
\textit{on $\mathcal{N}$, where the unbounded linear sectorial
operators $H(u,v)$ are normal in~$X$, the fields
$S,\,S^{-1}:\mathcal{N}\rightarrow\mathcal{L}(X)$ and
$T_{0}:\mathcal{N}\rightarrow\mathcal{L}(X^{\alpha},X)$ are
regular, and the field
$T_{0}:\mathcal{N}\rightarrow\mathcal{L}(X^{\alpha})$ is bounded.
In this case, if}
$$
a_{k}^{\beta}=o(\xi_{k})\; \, \quad (k \to +\infty),
\eqno (2.5)
$$
\textit{then the dynamics of equation} (2.1) \textit{is
finite-dimensional on the attractor.}
\medskip

\bigskip
\centerline{\large\bf 3. Parabolic systems}
\medskip

Now we consider the system of equations (1.1) on
$J=\mathbb{R}\,|\,\mathrm{mod}\,\mathbb{Z}$ with $u=(u_{1}
,...,u_{m})$. We assume that the matrix function $f=f(x,u)$ and
the vector function $g=g(x,u)$ belong to the smoothness class
$C^{\infty}$ on $J\times{\mathbb R}^{m}$ and write system (1.1) in
the abstract form (2.1) with $X=L^{2}(J,\mathbb{R}^{m})$, positive
definite operator $Au=u-Du_{xx}$, and nonlinearity $F:
u\rightarrow f(x,u)\partial_{x}u+g(x,u)$. Assume that
$\{X^{\alpha}\}_{\alpha\geq 0}$ is the Hilbert semiscale generated
by~$A$ and $\mathcal{H}^{s}=\mathcal{H}^{s}(J)$ are generalized
Sobolev $L^{2}$-spaces (spaces of Bessel potentials~[5,17]) of
scalar functions on~$J$ with arbitrary $s\geq 0$. If $s>1/2$, then
$\mathcal{H}^{s}\subset C(J)$ and $\mathcal{H}^{s}$ is a Banach
algebra [17, Section~2.8.3]. The differentiation operator
$\partial_{x}$ belongs to
$\mathcal{L}(\mathcal{H}^{s+1},\mathcal{H}^{s})$. As the phase
space we choose
$X^{\alpha}=\mathcal{H}^{2\alpha}(J,\mathbb{R}^{m})$ with
arbitrary $\alpha \in (3/4,1)$ which is fixed below.

We shall generalize the conclusions of [15, pp.~991--992] about
the smoothness of the nonlinear function~$F$ and the phase
dynamics of (1.1) to the case $m>1$. We let the symbol
$\hookrightarrow$ denote linear continuous embeddings of function
spaces and shall use necessary results obtained in~[5,17]. For an
arbitrary $C^{\infty}$-function
$z:J\times\mathbb{R}^{m}\rightarrow \mathbb{R}$, the mapping
$\psi:u\rightarrow z(x,u)$ is a function of class $BC^{\nu}$ from
$C^{s}(J)$ in $C^{s}(J)$ for all $\nu,s\in\mathbb{N}$. Since
$\mathcal{H}^{2\alpha}\hookrightarrow C^{1}(J)$, we have $\psi\in
BC^{\nu}(\mathcal{H}^{2\alpha},C^{1}(J))$. Embedding theorems
imply that $\psi\in
BC^{\nu}(\mathcal{H}^{s}(J,\mathbb{R}^{m}),\mathcal{H}^{s}(J))$.
As we see, $\,F\in BC^{2}(X^{1},X^{1/2})$ and $F\in
BC^{1}(X^{3/2},X^{1})$. Moreover, $X^{\alpha}\hookrightarrow
C^{1}(J,\mathbb{R}^{m})\hookrightarrow
C(J,\mathbb{R}^{m})\hookrightarrow X$, and hence $F\in
BC^{3}(X^{\alpha},X)$. We also note that $X^{3/2}\hookrightarrow
C^{2}(J,\mathbb{R}^{m})$ and $X^{2}\hookrightarrow
C^{3}(J,\mathbb{R}^{m})$.

In the case of finite functions $f=f(u),\,g=g(u)$, the
dissipativity of system (1.1) with phase space $X^{1/2}$, and
hence also with $X^{\alpha}$, $3/4<\alpha<1$, was proved in~[8,
Theorem~3.1]. This result can easily be transferred to the case of
functions $f(x,u)$ and $g(x,u)$ that are finite in~$u$ and can
also be generalized in other directions. Anyway, we further assume
that system (1.1) is dissipative in~$X^{\alpha}$ and there exists
the global attractor $\mathcal{A}\subset X^{\alpha}$. Using the
above-listed properties of nonlinearity~$F$ and following the
reasoning in~[15, p.~992], we formulate the following remark.
\medskip

\textbf{Remark~3.1} (see [15, Remark 5.2]). The following
assertions hold: (a) the attractor $\mathcal{A}$ is bounded
in~$X^{2}$; (b) if $Y$ is a Banach space, then each vector field
$\Pi:\mathcal{N}\rightarrow Y$ continuous in the
$(X^{\alpha}\times X^{\alpha})$-metric can be continued to
$C^{1}$-mapping $X^{1}\times X^{1}\rightarrow Y$ regularly in the
sense of Definition~2.1.

Our goal is to apply Theorem~2.2 to system~(1.1)
and to prove that the final dynamics is finite-dimensional.
Let
$$
G(u)=-Au+F(u)=D\partial_{xx}u-u +f(x,u)\partial_{x}u+g(x,u) \eqno (3.1)
$$
be the vector field of system (1.1), and let
$\mathcal{N}=\mathcal{A}\times\mathcal{A}\subset X^{\alpha}\times
X^{\alpha}$. The main idea, as in~[15], is related to the change
of variable in the linear differential expression with respect to
$x\in J$ for the difference $G(u)-G(v)$ for a fixed $(u,v)\in
\mathcal{N}$, which allows one to eliminate the dependence on
$\partial_{x}h$, $h(x)=u(x)-v(x)$. Along with the convection
matrix $f=\{f_{ij}\}$ we consider the $m\times m$ function
matrices
$$
g_{u}=\{\dfrac {\partial g_{i}}{\partial u_{j}}\},\quad \quad
f_{u}\partial_{x}u=\{\mathop{\sum}\limits_{l=1}\limits^{m}\dfrac
{\partial f_{il}}{\partial u_{j}}\partial_{x}u_{l}\},\quad \quad
i,j\in\overline{1,m}.
$$
We put
$$
B_{0}(x;u,v)=-E+\int_{0}^{1}(f_{u}(x,w(x))w_{x}(x)+g_{u}(x,w(x))d\tau,
\eqno (3.2.1)
$$
where $E$ is the unit $m\times m$ matrix, and
$$
B(x;u,v)=\int_{0}^{1}f(x,w(x))d\tau \eqno (3.2.2)
$$
for $u,v\in X^{\alpha}$, $w(x)=\tau u(x)+(1-\tau)v(x)$, $x\in J$.
The elements of the matrices $B_{0}$ and~$B$ are continuous
functions, and for $u,v\in\mathcal{A}$, function of class $C^{2}$
on~$J$.  If necessary, it is convenient to treat expressions (3.2)
as Bochner integrals ranging in some function spaces. Using the
$C^{1}$-smoothness of the mappings $(u,v)\rightarrow
f_{u}(x,w)w_{x}+g_{u}(x,w)$ and $(u,v)\rightarrow
f(x,w),\;X^{\alpha}\times X^{\alpha}\rightarrow
C(J,\mathbb{M}^{m})$ for a fixed $\tau\in[0,1]$ and
differentiating the expression under the integral sign in (3.2)
with respect to the parameter $(u,v)$, we conclude that the
mappings $(u,v)\rightarrow B_{0}(\cdot\,;u,v)$ and
$(u,v)\rightarrow B(\cdot\,;u,v)$ are of class
$C^{1}(X^{\alpha}\times X^{\alpha},C(J,\mathbb{M}^{m}))$. By the
integral mean-value theorem for nonlinear operators, we have
$$
G(u)-G(v)=-Ah+(\int_{0}^{1}\mathrm{D}F(\tau u+(1-\tau)v)d\tau)h
$$
$$
=Dh_{xx}+B_{0}(x;u,v)h+B(x;u,v)h_{x}\doteq Rh,
$$
where $h=u-v$, $u,v\in \mathcal{A}$, and $\tau u+(1-\tau)v\in {\rm
co}\,\mathcal{A}$. Here $\mathrm{D}$ is the Frechet
differentiation. Proceeding as in~[6], we apply the transformation
$h=U\eta$ to the differential expression $Rh$, where the $m\times
m$ matrix function $U(x)=U(x;u,v)$, $x\in [0,1]$, is a solution of
the linear Cauchy problem
$$
U_{x}=-\frac{1}{2} D^{-1} B(x)U,\quad \; U(0)=E.
\eqno (3.3)
$$
Similar problems are considered in [2, Ch. 3, 5]. We often write
$B_{0}$, $B$, and $U$ omitting the dependence on~$u$ and~$v$ and
sometimes on~$x$. Taking into account the fact that
$$
U_{xx}=-\frac{1}{2} D^{-1} (B_{x}(x)U+B(x)U_{x})
= -\frac{1}{2} D^{-1} B_{x} (x)U+\frac{1}{4} D^{-1} B(x)D^{-1} B(x)U,
$$
we have
$$
Rh=RU\eta =D(U\eta _{xx} +2U_{x} \eta _{x} +U_{xx} \eta )+B_{0}
(x)U\eta +B(x)(U_{x} \eta +U\eta _{x} )
$$
$$
=DU\eta _{xx}-B(x)U\eta _{x} -\frac{1}{2} B_{x} (x)U\eta +\frac{1}{4}
B(x)D^{-1} B(x)U\eta+B_{0} (x)U\eta +B(x)U\eta _{x}
$$
$$
+B(x)(-\frac{1}{2} D^{-1} B(x)U\eta ) =DU\eta _{xx} +(B_{0}
(x)-\frac{1}{2} B_{x} (x)-\frac{1}{4} B(x)D^{-1} B(x))U\eta.
$$
Now we write a decomposition of the form (2.2) for the vector field (3.1)
of evolution equation (1.1) on the attractor $\mathcal{A}$ with linear components
$$
T_{0}(u,v)h=\omega h+(B_{0}(x)-\frac{1}{2} B_{x} (x)-\frac{1}{4} B(x)D^{-1} B(x))h,
\eqno (3.4.1)
$$
$$
T(u,v)h=\omega h-DU\partial _{xx} U^{-1} h,
\eqno (3.4.2)
$$
where the numerical parameter $\omega >0$ will be chosen later.

Everywhere below, $I\doteq {\rm Id}$ in a Banach space. By
$\mathbb{M}^{m}$ we denote the algebra of numerical $m\times m$
matrices with Euclidean norm, and by $Y(J,\mathbb{M}^{m})$ we
denote the linear spaces of such matrices with elements from some
Banach space~$Y$ of scalar functions on~$J$  or $[0,1]$ . We
slightly generalize the fact that linear problem (3.3) can be
solved explicitly under the condition that the operators
$D^{-1}B(x)$ are commutative in $x\in J$.
\medskip

\textbf{Lemma 3.2.}
\textit{Let $D^{-1}B(x)=CW(x)C^{-1}$
with constant nondegenerate matrix~$C$ and matrix function
$W\in C(J,\mathbb{M}^{m})$,
and let $W(x_{1})W(x_{2})=W(x_{2})W(x_{1})$ for $x_{1},x_{2}\in J$.
Then $U(x)=C\overline{U}(x)C^{-1}$ with}
$$
\overline{U}(x)=\mathrm{exp}(-\frac{1}{2}\int\limits_0^{x}W(\xi)d\xi\,)
\eqno (3.5)
$$
\textit{is a solution of the Cauchy problem $(3.3)$ on $[0,1]$} .
\medskip

\textbf{Proof.}
Under the conditions of the lemma, we have
$\overline{U}(x)W(x)=W(x)\overline{U}(x)$,
and hence $\overline{U}_{x}=-\frac{1}{2}W(x)\overline{U}$.
Further, $\overline{U}(0)=U(0)=E$ and
$$
U_{x}=C\overline{U}_{x}C^{-1}=C(-\frac {1} {2}W\overline{U})C^{-1}
$$
$$
=C(-\frac {1} {2}C^{-1}D^{-1}BC\cdot C^{-1}UC)C^{-1}=-\frac {1} {2}D^{-1}BU.
\quad \Box
$$

Now we prove regularity in the sense of Definition~2.1 of some vector fields
on the compact set $\mathcal{N}\subset X^{\alpha}\times X^{\alpha}$.
If $Y\hookrightarrow Y_{1}$ for the function spaces $Y$ and $Y_{1}$,
then the regularity of the field $\Pi:\mathcal{N}\rightarrow Y$
implies the regularity of $\Pi:\mathcal{N}\rightarrow Y_{1}$.
\medskip

\textbf{Lemma 3.3.}
\textit{The field of operators $T_{0}$ on $\mathcal{N}$
is bounded ranging in $\mathcal{L}(X^{\alpha})$
and regularly ranging in $\mathcal{L}(X^{\alpha},X)$.}
\medskip

\textbf{Proof.}
Let $\|\cdot\|_{\alpha,\alpha}$ and $\|\cdot\|_{\alpha,0}$
be norms of operators in $\mathcal{L}(X^{\alpha})$
and $\mathcal{L}(X^{\alpha},X)$.
We assume that $T_{0}h=Q(x;u,v)h$ in (3.4.1)
with $h\in\mathrm{co}\,\mathcal{A} \subset X^{\alpha}$.
By Remark~3.1.(a), the convex hull of the attractor
is bounded in the norm of~$X^{2}$ which is equivalent to the norm
$\mathcal{H}^{4}(J,\mathbb{R}^{m})$,
and hence the matrix functions $B$, $B_{0}$, and $BD^{-1}B$
are uniformly bounded with respect to $(u,v)\in \mathcal{N}$
in $\mathcal{H}^{3}(J,\mathbb{M}^{m})$.
Thus, the matrix functions $B_{x}$ and $Q$ are bounded on $\mathcal{N}$
in the norm $\mathcal{H}^{2}(J,\mathbb{M}^{m})$
and $T_{0}$ is the operator of multiplication of vector functions in
$X^{\alpha}=\mathcal{H}^{2\alpha}(J,\mathbb{R}^{m})$
by the matrix $Q\in \mathcal{H}^{2\alpha}(J,\mathbb{M}^{m})$
with  $2\alpha\in(3/2,2)$.
Since $\mathcal{H}^{2\alpha}(J)$ is a Banach algebra,
we see that $T_{0}(u,v)\in \mathcal{L}(X^{\alpha})$ and
$\|T_{0}(u,v)\|_{\alpha,\alpha}\leq \mathrm{const}$ on~$\mathcal{N}$.

Since $\mathcal{H}^{2\alpha}(J)\hookrightarrow C(J)\hookrightarrow
L^{2}(J)$, we have
$\mathcal{H}^{2\alpha}(J,\mathbb{M}^{m})\hookrightarrow
C(J,\mathbb{M}^{m})\hookrightarrow L^{2}(J,\mathbb{M}^{m})$ and
$\|T_{0}\|_{\alpha,0}\leq c\|Q\|_{0,0}$, where $\|Q\|_{0,0}$ is
the norm of~$Q$ as an operator in $\mathcal{L}(X)$ and
$c=c(\mathcal{A})$. Therefore, the field of operators
$T_{0}:\mathcal{N}\rightarrow\mathcal{L}(X^{\alpha},X)$ is regular
if the field of matrix functions $Q:\mathcal{N}\rightarrow
L^{2}(J,\mathbb{M}^{m})$ is regular. The function $u\rightarrow
f(x,u),\;X^{\alpha}\rightarrow C(J,\mathbb{M}^{m})$, is of class
$C^{1}$. Since the mappings $(u,v)\rightarrow B_{0}(\cdot\,;u,v)$
and $(u,v)\rightarrow B(\cdot\,;u,v)$ are of class
$C^{1}(X^{\alpha}\times X^{\alpha},C(J,\mathbb{M}^{m}))$, it
follows that their restrictions to $\mathcal{N}$ are regular. The
regularity of the field $BD^{-1}B:\mathcal{N}\rightarrow
C(J,\mathbb{M}^{m})$ follows from the regularity of the fields~$B$
and $D^{-1}=\mathrm{const}$ with the multiplicative structure of
$C(J,\mathbb{M}^{m})$ taken into account. Moreover, the fields of
matrix functions $B$, $B_{0}$, $BD^{-1}B$ on $\mathcal{N}$ are
regular with values in $L^{2}(J,\mathbb{M}^{m})$.

Now we prove the regularity of the field $\Pi\doteq
B_{x}:\mathcal{N}\rightarrow L^{2}(J,\mathbb{M}^{m})$. Let
$\Pi_{\tau}(u,v)=(f(x,w))_{x}$ with $w=\tau u(x)+(1-\tau)v(x)$ for
a fixed $\tau\in [0,1]$ and arbitrary $u=u(x),\,v=v(x)\in X^{1}$.
Then $\Pi(u,v)=(B(x;u(x),v(x)))_{x}$ is the result of integration
of $\Pi_{\tau}(u,v)$ over $\tau$. The mapping $u\rightarrow
f(x,u)$ belongs at least to the class $BC^{1}(X^{1},X^{1/2})$, and
hence $\Pi_{\tau}\in C^{1}(X^{1}\times
X^{1},L^{2}(J,\mathbb{M}^{m}))$. Differentiating the integral
expression for $\Pi(u,v)$ with respect to the parameter $(u,v)\in
X^{1}\times X^{1}$, we obtain $\Pi\in C^{1}(X^{1}\times
X^{1},L^{2}(J,\mathbb{M}^{m}))$. It remains to verify that the
fields $\Pi:\mathcal{N}\rightarrow L^{2}(J,\mathbb{M}^{m})$ are
continuous and to use Remark~3.1.(b). By~[15, Lemma~1.1], the
function $u\rightarrow Au$, $\mathcal{A}\rightarrow X$, with
$Au=u-Du_{xx}$ is continuous in the $X^{\alpha}$-metric; the same
holds for the mappings $u\rightarrow u_{xx}$ and $u\rightarrow
u_{x}$ of the set $\mathrm{co}\,\mathcal{A}\subset X^{\alpha}$
into $X$ for $u\in\mathrm{co}\,\mathcal{A}\subset X^{1}$. In the
relation $(f(x,u))_{x}=f_{x}+f_{u}u_{x}$, the operators
$u\rightarrow f_{x}(x,u)$, $u\rightarrow f_{u}(x,u)u_{x}$
continuously act from $X^{\alpha}$ to $C(J,\mathbb{M}^{m})$, and
hence $\Pi_{\tau},\,\Pi \in
C(\mathcal{N},L^{2}(J,\mathbb{M}^{m}))$. The proof of the lemma is
complete. $\Box$

Everywhere below, $I\doteq {\rm Id}$ in a Banach space. The matrix
functions $B(x)$ and $U(x)$ in the Cauchy problem (3.3) can be
treated as bounded linear operators in~$X$. The following
assertion is related to the smooth dependence of solutions of
differential equations on a parameter.
\medskip

\textbf{Lemma 3.4.} \textit{The field of operators
$U:\mathcal{N}\rightarrow \mathcal{L}(X)$ is regular.}
\medskip

\textbf{Proof.}
We consider (3.3) for arbitrary $u,v\in X^{\alpha}$
as the non-autonomous evolution problem
$$
\partial_{x}U=-\frac{1}{2} D^{-1}B(x;(u,v))U,\quad \; U(0)=I
$$
in the Banach algebra $\mathcal{L}(X)$ with identically zero sectorial linear part
and the parameter $(u,v)\in X^{\alpha}\times X^{\alpha}$.
The function
$$
(x,U,(u,v))\rightarrow -\dfrac{1}{2}D^{-1}B(x;(u,v))U
$$
ranging in $\mathcal{L}(X)$ is Lipschitz in $x$, linear in $U\in
\mathcal{L}(X)$ and of class~$C^{1}$ with respect to the parameter
$(u,v)$. Under these conditions, by [5, Theorem~3.4.4], the
mapping $(u,v)\rightarrow U(x;(u,v)),\;X^{\alpha}\times
X^{\alpha}\rightarrow \mathcal{L}(X)$ is continuously
differentiable, and hence the operator field
$U:\mathcal{N}\rightarrow \mathcal{L}(X)$ is regular. $\Box$

Now we formulate an important condition on the diffusion matrix
$D$ and the convection matrix $f$ of system (1.1).
\medskip

\textbf{Assumption~3.5.} \textit{$Df(x,u)=f(x,u)D$ for $x\in
J,\,u\in {\rm co}\,\mathcal{A}$.}
\medskip

For the scalar diffusion matrix $D=dE$, this assumption is
satisfied automatically. In the case of $m$ distinct diffusion
coefficients $d_{j}$, Assumption~3.5 holds under the condition
that the matrix $f$ is diagonal on $\mathrm{co}\,\mathcal{A}$, and
in the case of~$s$ distinct diffusion coefficients, $1<s<m$, it
holds under the condition that the matrix $f$ on
$\mathrm{co}\,\mathcal{A}$ inherits the block structure (with
respect to the same $d_{j}$) of the matrix
$D=\mathrm{diag}\{d_{j}\}$.
\medskip

{\bf  Lemma 3.6.} \textit{If Assumption~$3.5$ holds, then}
$$
T(u,v)=U(u,v)(\omega I-D\partial_{xx})U^{-1}(u,v)
$$
\textit{for $u,v \in \mathcal{A}$.}
\medskip

\textbf{Proof.} Assumption~3.5 implies (for any $x\in J$ and
$u,v\in \mathcal{A}$) that $DB(x)=B(x)D$ for the matrices
$B(x)=B(x;u,v)$ in~(3.2.2). Thus, the matrices $B(x)$ and
$D^{-1}B(x)$ inherit the block structure (with respect to the same
$d_{j}$) of the diffusion matrix
$D=\mathrm{diag}\,\{d_{1},\dots,d_{m}\}$. Therefore, the same also
holds for the solutions $U(x)$ of problem (3.3), and hence
$DU(x)=U(x)D,\;x\in [0,1]$, and the assertion of the lemma follows
from (3.4.2). $\Box$

\bigskip
\centerline{\large\bf 4. Main results}
\medskip

The conditions for the dynamics to be finite-dimensional on the
attractor will depend on the structure of the diffusion matrix~$D$
and the nonlinear function~$f$ in~(1.1). By Theorem~2.2, we need
to prove that the operators $T(u,v)$ in (3.4.2) are ``uniformly
and regularly'' similar, like (2.4), to the normal operators
in~$X$ and to establish the required sparseness (2.5) of the total
spectrum~$\Sigma_{T}$.

We note that $B(0)=B(1)$, $B_{x}(0)=B_{x}(1)$ for the matrix
function $B(x)=B(x;u,v)$ in (3.2.2) defined on $J\times
\mathcal{N}$. The matrix function $V(x)=U^{-1}(x),\;x\in [0,1]$,
is a solution (see~[2, Sect.~3.1.3]) of the Cauchy problem adjoint
to (3.3):
$$
V_{x} =\frac{1}{2} VD^{-1} B(x),\quad \; V(0)=E,
\eqno (4.1)
$$
and we have
$$
\eta=Vh, \quad \eta_{x} =V_{x}h+Vh_{x},\quad V_{x}=\frac{1}{2}VD^{-1}B,
$$
$$
V_{x}(0)=\frac{1}{2}V(0)D^{-1}B(0), \; \; V_{x}(1)=\frac{1}{2} V(1)D^{-1}B(1),
$$
$$
\eta(0)=h(0),\quad \eta (1)=V(1)h(1),
$$
$$
\eta_{x}(0)=\frac{1}{2}D^{-1} B(0)h(0)+h_{x} (0),\quad \eta_{x} (1)=\frac{1}{2}
V(1)D^{-1} B(1)h(1)+V(1)h_{x}(1).
$$
So the periodic boundary conditions $h(1)=h(0),\;h_{x}(1)=h_{x}(0)$ become
$$
\eta (1)=V(1)\eta(0),\quad \quad \quad \eta _{x} (1)=V(1)\eta_{x}(0),
\eqno (4.2)
$$
where $V(1)\neq E$ in general. Further, we use the notation
$U_{1}=U_{1}(u,v),\;V_{1}=V_{1}(u,v)$ with $(u,v)\in\mathcal{N}$
for the monodromy operators $U(1),\;V(1)\in \mathcal{L}({\mathbb
R}^{m})=\mathbb{M}^{m}$.

The following assertion plays the key role.
\medskip

\textbf{Lemma 4.1.} \textit{If system} (1.1) \textit{is
dissipative in $X^{\alpha}$ with $\alpha \in (3/4, 1)$, then the
phase dynamics on the attractor is finite-dimensional in each of
the following two cases}.

(i) \textit{The diffusion matrix $D$ is scalar and, for all
$u,v\in\mathcal{A}$, the monodromy operators $V_{1}(u,v)$ are
similar to positive definite ones with a fixed similarity matrix
$C=C(\mathcal{A})$}.

(ii) \textit{Assumption $3.5$ holds and, for all
$u,v\in\mathcal{A}$, the monodromy operators $V_{1}(u,v)$ are
similar to diagonal positive definite ones with a fixed similarity
matrix $C=C(\mathcal{A})$}.
\medskip

\textbf{Proof.} By the conditions of the lemma, we have
$V_{1}=C^{-1}\mathcal{V}C$ for positive definite operators
$\mathcal{V}=\mathcal{V}(u,v)$ in $\mathbb{R}^{m}$. For fixed
$u,v\in \mathcal{A}$, we let $\varphi_{j}\in{\mathbb R}^{m}$ and
$\mu_{j}>0$ denote orthonormal eigenvectors and eigenvalues of the
operator $\mathcal{V}$ with $j\in\overline{1,m}$. We assume that
$H_{0}=H_{0}(u,v)=\omega
I-D\partial_{xx},\;D=\mathrm{diag}\{d_{j}\}$, with boundary
conditions (4.2) on $(0,1)$ for some $\omega>0$. We also assume
that
$$
\chi_{k,j}(x)=e^{2\pi kix} \cdot \varphi _{j},\quad x\in J,\quad k\in {\mathbb Z},\quad j\in\overline{1,m}.
$$
Since $\mathcal{V} \varphi_{j}=\mu_{j}\varphi_{j}$,
we have $V_{1}C^{-1}\varphi_{j}=\mu_{j}C^{-1}\varphi_{j}$
and, for the functions $\psi_{k,j}(x)=\mu_{j}^{x}\cdot C^{-1}\chi_{k,j}$,
$\psi_{k,j}(0)=C^{-1}\varphi_{j}$, $\psi_{k,j}(1)=V_{1}C^{-1}\varphi_{j}$,
$$
(\psi_{k,j})_{x}(0)=(\ln\mu_{j}+2\pi
ki)C^{-1}\varphi_{j},\quad(\psi_{k,j})_{x}(1) =(\ln\mu_{j}+2\pi
ki)V_{1}C^{-1}\varphi_{j}.
$$
As we see, $\psi_{k,j}$ are eigenfunctions of the operator $H_{0}$ with eigenvalues
$$
\lambda_{k,j} =\omega-d_{j}(\ln \mu_{j}+2\pi ki)^{2} =\omega
+d_{j} (2\pi k-i\ln\mu_{j})^{2}, \eqno (4.3)
$$
where $d_{j}\equiv d>0$ in case (i). The operators $V_{1}(u,v)$
continuously depend on $(u,v)\in \mathcal{N}$, and hence this also
holds for their spectrum. By the compactness of
$\mathcal{N}\subset X^{\alpha}\times X^{\alpha}$, we have $0<c_{1}
\le \mu _{j} \le c_{2}$, $j\in\overline{1,m}$, for some
$c_{1}({\mathcal{A}}),\,c_{2}({\mathcal{A}})$. Thus, the values
$\left|\,\ln\mu _{j}\right|$ are uniformly bounded in
$j\in\overline{1,m}$ and $u,v\in\mathcal{A}$. We put
$$
S_{0}(x)=CV_{1}^{-x}=\mathcal{V}^{-x}C
$$
for $x\in J$ and $H=S_{0}H_{0}S_{0}^{-1}$.
Then $S_{0}\psi_{k,j}=\chi_{k,j}$ and
$$
H\chi_{k,j}=S_{0}H_{0}\psi_{k,j}=\lambda_{k,j}\,S_{0}\psi_{k,j}=\lambda_{k,j}\,\chi_{k,j}.
$$
Since the system of functions $\{\chi_{k,j}\}$ is complete and orthonormal in
$X=L^{2}(J,{\mathbb R}^{m})$, it follows that the operators $H=H(u,v)$
are normal in~$X$ for $u,v\in \mathcal{A}$.
Let $S=S_{0}U^{-1}(x)=S_{0}V(x)$.
We use Lemma~3.6 to write decomposition~(2.4) of the vector field (3.1)
on the attractor $\mathcal{A}$ with
$$
T(u,v)=UH_{0}U^{-1}=US_{0}^{-1}HS_{0}U^{-1}=S^{-1}(u,v)H(u,v)S(u,v)
$$
and operators $T_{0}(u,v)$ of the form~(3.4.1).
We see that $S^{-1}=U(x)S_{0}^{-1}$.

By Lemma 3.4, the operator field $U$ on $\mathcal{N}$ ranging in
the Banach algebra $\mathcal{L}(X)$ is regular, and hence, the
field of inverse operators
$V:\mathcal{N}\rightarrow\mathcal{L}(X)$ is regular. Since
$V_{1}=V(1)$ and $\mathcal{V}=CV_{1}C^{-1}$, it follows that the
operator field $\mathcal{V}:\mathcal{N}\rightarrow\mathbb{M}^{m}$
is regular and $\|\partial_{t}\mathcal{V}(u,v)\|\leq c_{3}$ for
the derivative $\partial_{t}\mathcal{V}$ at zero for all $(u,v)\in
\mathcal{N}$ (see Definition~2.1). Here $\|\cdot\|$ is the
Euclidean norm of matrices. Let $b=2\max (c_{2},c_{3})$ and
$\delta=1-c_{1}/b$, then $\delta\in (0,1)$. Since the spectrum
$\sigma(b^{-1}\mathcal{V}-E)\subset (-\delta,0)$ and
$\|b^{-1}\mathcal{V}-E\|<\delta$, it follows that the matrix
representation
$$
\ln \mathcal{V}=\ln(bE)+\ln (E+b^{-1} \mathcal{V}-E)=\ln(bE)
+\mathop{\sum}\limits_{n=1}\limits^{\infty}\dfrac{(-1)^{n-1}}{n}(b^{-1}
\mathcal{V}-E)^{n} \eqno(4.4)
$$
converges uniformly on $\mathcal{N}$. By~[11, Sect.~5.8,
Exercise~3], we have
$$
\partial_{t}(b^{-1}
\mathcal{V}-E)^{n}=\mathop{\sum}\limits_{i=1}\limits^{n}(b^{-1}
\mathcal{V}-E)^{i-1}\partial_{t}(b^{-1}
\mathcal{V})(b^{-1}\mathcal{V}-E)^{n-i},
$$
and therefore, $\|\partial_{t}(b^{-1}\mathcal{V}-E)^{n}\|<
n\delta^{n-1}$. If we differentiate~(4.4) with respect to~$t$, we
obtain a uniformly converging series with the estimate
$$
\|\partial_{t}\ln \mathcal{V}(u,v)\|<
\mathop{\sum}\limits_{n=1}\limits^{\infty}\delta^{n-1}=\dfrac{1}{1-\delta}\,,\quad
(u,v)\in \mathcal{N}.
$$
So the operator field $\ln
\mathcal{V}:\mathcal{N}\rightarrow\mathbb{M}^{m}$ is regular. We
have $\mathcal{V}^{-x}=\mathrm{exp}(-x\ln \mathcal{V})$ and the
standard decomposition of the matrix exponent guarantees that the
field $\mathcal{V}^{-x}:\mathcal{N}\rightarrow \mathcal{L}(X)$ is
regular, which implies the regularity of fields of the operators
$S_{0},S,S^{-1}:\mathcal{N}\rightarrow\mathcal{L}(X)$.

Finally, by Lemma~3.3, the field of operators $T_{0}$ in (3.4.1)
is regular on $\mathcal{N}$ with values in
$\mathcal{L}(X^{\alpha},X)$ and bounded with values in
$\mathcal{L}(X^{\alpha})$.

Let $\Sigma_{H}=\Sigma_{T}$ be the total spectrum of the field of
operators $H(u,v)$ on $\mathcal{N}$. Using (4.3), we choose a
parameter $\omega>0$ that ensures the inclusion $\Sigma_{H}
\subset \Omega (c,\theta)$ of the form (2.3) with $\theta =1/2$
and an appropriate $c>0$. As we can see, operators $H(u,v)$ are
sectorial. Since $3/4<\alpha <1$, we have $\beta=(\alpha +\theta
)/3<1/2$. Moreover, from (4.3) we derive that the set ${\mathbb
C}\backslash\Sigma_{H}$ contains vertical strips
$\Gamma(a_{k},\xi_{k})$ with
$$
a_{k} \sim 4\pi ^{2} k^{2} ,\; \; \; \xi _{k} \sim 4\pi ^{2}k
$$
and hence $a_{k}^{\beta}=o(\xi_{k})$ as $k\rightarrow\infty$.
Thus, all conditions of Theorem~2.2 are satisfied and the dynamics
of system (1.1) is finite-dimensional on the attractor. $\Box$
\medskip

\textbf{Remark 4.2.}
For all $u,v\in \mathcal{A}$, the monodromy operators $U_{1}=U_{1}(u,v)$
and $V_{1}=V_{1}(u,v)=U_{1}^{-1}$ are positive definite
if any of the following two conditions is satisfied:

(a) the matrices $D^{-1}B(x_{1})$ and $D^{-1}B(x_{2})$
are symmetric and commutative for $x_{1},x_{2}\in J$;

(b) $(D^{-1}B(x))^{\mathrm{t}}=D^{-1}B(1-x)$ for all $x\in J$,
where $(\cdot)^{\mathrm{t}}$ is the operation of transposition.

Under conditions~(a), Lemma~3.2 holds with $C=E$ and the matrix
$$
U_{1}=\mathrm{exp}(-\frac {1} {2}\int\limits_0^{1}D^{-1}B(x)dx\,)
$$
is positive definite.
The sufficiency of condition~(b) was proved in~[18, Proposition~2.3].
\medskip

\textbf{Theorem 4.3.} \textit{Assume that system $(1.1)$ is
dissipative in $X^{\alpha}$ with $\alpha \in (3/4,\,1)$ and the
convection matrix $f$ is diagonal on $\mathrm{co}\,\mathcal{A}$.
Then the phase dynamics is finite-dimensional on the attractor.}
\medskip

\textbf{Proof}. Under the conditions of the theorem, the matrices
$B(\cdot; u, v)$ from (3.2.2), and hence (see the proof of Lemma
3.6), also the matrices $U(\cdot;u,v),V(\cdot;u,v)$, are diagonal
on $\mathcal{A}\times \mathcal{A}$. According to Remark 4.2.(a)
matrices $U(\cdot; u,v),V (\cdot;u,v)$ are positive definite. The
monodromy operators $V_{1}(u,v)$ are also positive definite and
are diagonal, so we can refer to Lemma 4.2.(ii) with $C=E$.

\medskip
\textbf{Lemma 4.4.} \textit{Assume that system $(1.1)$ is
dissipative in $X^{\alpha}$ with $\alpha \in (3/4,\,1)$ and
$D=dE$. Then the phase dynamics is finite-dimensional on the
attractor if, for $(x,u)\in {J\times\rm co\,\mathcal{A}}$},
$$
D^{-1}f(x,u(x))=CH(x,u(x))C^{-1}, \eqno (4.5)
$$
\textit{where the symmetric matrix functions $H(x;u)$ commute with
each other for any $(x,u)\in J\times \mathrm{co}\,\mathcal{A}$ and
$C$ is a constant nondegenerate matrix.}
\medskip

\textbf{Proof.} From (3.2.2) we derive that
$D^{-1}B(x)=CW(x)C^{-1}$, where
$$
W(x)=W(x;u,v)=\int_{0}^{1}H(x;w(x))d\tau
$$
for $u,v\in \mathcal{A}$, $w(x)=\tau u(x)+(1-\tau)v(x)$, $x\in J$.
By Lemma~3.2, the monodromy operator $V_{1}(u,v)=U_{1}^{-1}(u,v)$,
$u,v\in \mathcal{A}$, satisfies the relation
$V_{1}=C(\overline{U}(1))^{-1}C^{-1}$ with operator $\overline{U}$
given in formula~(3.5). In this case,
$\overline{U}_{x}=-\dfrac{1}{2}W(x)\overline{U}$ and the matrices
$W(x;u,v)$ are symmetric and commutative on~$J$ for all $u,v\in
\mathcal{A}$. By Remark~4.2.(a), the operator $\overline{U}(1)$ is
positive definite and the assertion of the theorem follows from
Lemma~4.1.(i). $\Box$

We shall give two more arguments ensuring that the final dynamics
is finite-dimensional.

\medskip
\textbf{Theorem 4.5.} \textit{Assume that system $(1.1)$ is
dissipative in $X^{\alpha}$ with $\alpha \in (3/4,\,1)$ and
$D=dE$. Then the phase dynamics is finite-dimensional on the
attractor if the following two conditions are satisfied}:

(i) \textit{the numerical matrices $f(x,u(x))$ have $m$ distinct
real eigenvalues for each $(x,u)\in J\times
\mathrm{co}\,\mathcal{A}$};

(ii) \textit{the matrices $f(x,u)$ commute with each other for any
$(x,u)\in J\times \mathrm{co}\,\mathcal{A}$}.
\medskip

\textbf{Proof.} Condition (ii) and assumption $D=dE$ imply that
the matrices $D^{-1}f(x,u(x))$ commute with each other on
$J\times\mathrm{co}\,{\mathcal{A}}$. It is known [11,
Theorem~8.6.1] that two simple (similar to diagonal) commutative
$m\times m$ matrices have a common set of~$m$ of linearly
independent eigenvectors. By condition~(i), all eigenvalues of
each numerical matrix $f(x,u(x))$ with $(x,u)\in J\times
\mathrm{co}\,\mathcal{A}$ are real and distinct, and hence there
exists a unique (up to permutations and multiplications by~$-1$)
common (for all these matrices) normalized basis
$\mathcal{E}=(e_{1},\dots,e_{m})$ of their eigenvectors
in~$\mathbb{R}^{m}$. By~$C$ we denote the constant matrix of
transition from the canonical basis in~$\mathbb{R}^{m}$ to the
basis $\mathcal{E}$, and by~$H(x)$ we denote diagonal (symmetric)
matrices of linear operators $D^{-1}f(x,u(x))\in
\mathcal{L}(\mathbb{R}^{m})$ in this basis. We see that relation
(4.5) is satisfied and it remains to apply Lemma~4.4. $\Box$
\medskip

\textbf{Theorem 4.6.} \textit{Assume that system $(1.1)$ is
dissipative in $X^{\alpha}$ with $\alpha \in (3/4,\,1)$ and
$D=dE$. Then the phase dynamics is finite-dimensional on the
attractor if the matrices $f(x,u)$ are symmetric and commute with
each other for any $(x,u)\in J\times \mathrm{co}\,\mathcal{A}$.}
\medskip

\textbf{Proof.} The conditions of the theorem guarantee that the
numerical matrices

\noindent$D^{-1}f(x,u(x))$ commute with each other on $J\times
\mathrm{co}\,\mathcal{A}$. As in the proof of Lemma~4.4, from
formula (3.2.2) for $B(x)$, we derive that the matrices
$D^{-1}B(x_{1})$ and $D^{-1}B(x_{2})$ are symmetric and
commutative for arbitrary $x_{1},x_{2}\in J$. By Remark 4.2.(a),
the monodromy operators $V_{1}(u,v)$ are positive definite for any
$u,v\in \mathcal{A}$ and the assertion of the theorem follows from
Lemma~4.1.(i) with $C=E$.
\medskip

In contrast to Theorem~4.5, we here admit the multiplicity of
eigenvalues of the numerical matrices $f(x,u(x))$, but we assume
that these matrices are symmetric.

\bigskip
\centerline{\large\bf 5. Some examples}
\medskip

We consider several examples illustrating the above-described
theory in terms of properties of the convection matrix~$f$. Here
we restrict ourselves to the case of scalar diffusion and assume
that system (1.1) is dissipative in the phase space $X^{\alpha}$
with $\alpha\in (3/4,1)$. We assume that all the conditions
assumed below on $f=f(x,u)$ are valid for $x\in J$ and $u=u(x)$,
$u\in\mathrm{co}\,\mathcal{A}$.

\textbf{Proposition 5.1.} \textit{Assume that $D=dE$ and
$f(x,u)=f_{1}(x,u)Q$ with a scalar $C^{\infty}$-function $f_{1}$
and numerical $m\times m$ matrix~$Q$. Then, the dynamics on the
attractor of system $(1.1)$ is finite-dimensional if any of the
following two conditions is satisfied}:

(i) \textit{the matrix $Q$ has $m$ distinct real eigenvalues and
$f_{1}(x,u(x))\neq 0$ for $x\in J$ and
$u\in\mathrm{co}\,\mathcal{A}$};

(ii) \textit{the matrix $Q$ is symmetric.}
\medskip

\textbf{Proof.} The numerical matrices $f=f_{1}(x,u(x))Q$ are
commutative. In the case of~(i), each of these matrices has
distinct real eigenvalues $\lambda_{j}f_{1}(x,u(x))$, where
$\lambda_{1},\dots,\lambda_{m}$ are eigenvalues of~$Q$, and
Theorem~4.5 can be applied. In the case of~(ii), the fact that the
dynamics is finite-dimensional on the attractor is a direct
consequence of Theorem~4.6. $\Box$

\textbf{Remark 5.2.}
Condition~(i) in Proposition~5.1 is satisfied on~$Q$
for upper-triangular and lower-triangular matrices with distinct elements
on the diagonal.
For $m=2$ and $Q=\{q_{jl}\}$, this condition precisely means that
$(q_{11}-q_{22})^{2}+4q_{12}q_{21}>0$.
\medskip

\textbf{Example 5.3.} The dynamics on the attractor of system
(1.1) is finite-dimensional in the case of $m=2$, $D=dE$ and
$f(x,u)=\{f_{jl}(x,u)\}$ with $f_{11}=f_{22}$ and $f_{12}=f_{21}$.
This is a consequence of Theorem~4.6 and the commutativity of
numerical matrices of the form $\left(\begin{array}{cc} {a} & {b}
\\ {b} & {a} \end{array}\right)$.
\medskip

\textbf{Example 5.4.} Assume that $D=dE$, the matrix $f=P_{n}(Q)$,
where $P_{n}$ is a polynomial of degree $n\geq 0$ with
coefficients $a_{i}=a_{i}(x,u)$, $0\leq i \leq n$, $a_{i}\in
C^{\infty}(J\times \mathbb{R}^{m},\mathbb{R})$, and the numerical
matrix~$Q$ is symmetric. Then the dynamics of the attractor of
system (1.1) is finite-dimensional. This is also a consequence of
Theorem~4.6.
\medskip

\textbf{Proposition 5.5.} \textit{Assume that $D=dE$ and $f=Q(x)$,
where $Q$ is a $C^{\infty}$ function matrix. Then the dynamics of
system~$(1.1)$ is finite-dimensional on the attractor if
$Q^{\,\mathrm{t}}(x)=Q(1-x)$ for $x\in J$.}
\medskip

\textbf{Proof.} Since $f=Q(x)$, the matrix $B(x)$ in~(3.2.2)
satisfies the condition $(D^{-1}B(x))^{\mathrm{t}}=D^{-1}B(1-x)$
for all $x\in J$ and $u,v\in \mathcal{A}$. By Remark 4.2.(b), the
monodromy operators $V_{1}(u,v)$ are positive definite for all
$u,v\in \mathcal{A}$, and we can apply Lemma~4.1.(i) with $C=E$.
$\Box$
\medskip

\bigskip
\centerline{\large\bf 6. Other possible approaches}
\medskip

The above presentation is based on Theorem~2.2, which means the
verification of regularity (in the sense of Definition~2.1) of
operator vector fields on the attractor. Alternatively, one can
obtain a finite-dimensional reduction of one-dimensional parabolic
systems by using the technique~[19, Sect.~2.3, 2.4, 3.3] closely
related to the results obtained in~[4] about the dichotomies of
\textit{non-autonomous} parabolic equations. Here we will discuss
the $X^{1/2}$-dissipative systems of general form~(1.2). In our
short description (on the sketch level), we omit technical details
and refer to the criteria for the final dynamics to be
finite-dimensional, i.e., criteria~(Fl) and~(GrF) in Section~2.
The results of~[16, Sect.~3.6] about inverse uniqueness of
solutions of SPE~(2.1) allow one to conclude that the phase
semiflow on the attractor $\mathcal{A}$ expands to the continuous
flow $\{\Phi_{t}\}_{t\in \mathbb{R}}$. If
$u_{1},u_{2}\in\mathcal{A}$ and $h(t)=\Phi_{t}u_{1}-\Phi_{t}u_{2}$
for $t\in\mathbb{R}$, then
$$
h_{t}=Dh_{xx}+B_{0}(t,x)h+B(t,x)h_{x},
\eqno (6.1)
$$
$$
B_{0}(t,x)=\int_{0}^{1}f_{u} (x,w,w_{x})d\tau,\quad \;
B(t,x)=\int_{0}^{1}f_{u_{x}}(x,w,w_{x})d\tau,
$$
$w=\tau\Phi_{t}u_{1}+(1-\tau )\Phi_{t}u_{2}$, with matrix
functions $B_{0}$ and~$B$ sufficiently smooth in $(t,x)\in
\mathbb{R}\times J$ and bounded in $(u_{1},u_{2})$. We assume that
the invertible change $v(t,x)=S(t,x)h(t,x)$ with operators
$S,S^{-1}\in \mathcal{L}(X)$ depending on $u_{1},u_{2}$ allows one
to reduce~(6.1) to the equation
$$
v_{t}=-Hv+R(t)v,
$$
where $H=H(u_{1},u_{2})\in\mathcal{L}(X^{1},X)$ are normal (or
uniformly with respect to $(u_{1},u_{2})$ similar to normal)
sectorial operators and
$R=R(\cdot\,;u_{1},u_{2}):\mathbb{R}\rightarrow \mathcal{L}(X)$ is
a continuous operator function. Assume that the norms of the
operators $S$, $S^{-1}$, and $R$ are uniformly bounded in the
parameter $(u_{1},u_{2})$. In this case,  if the spectrum
$\Sigma_{H}$ combined over $u_{1},u_{2}\in \mathcal{A}$ is
``sufficiently rare'', then using the technique give in~[19], one
can verify that the phase flow is Lipschitzian on the attractor
and then apply criterion~(Fl). In contrast to the preceding
presentation, we here have to deal with second-order linear
differential expression in $t\in \mathbb{R}$ and not in $x\in J$.
The assertions of Section~4 can be obtained in this way after the
change $v(t,x)=V^{-x}(t,1)V(t,x)h(t,x)$, where $V(\cdot,x)$ is a
solution of the Cauchy problem (4.1).

Another possible approach to the problem of finite-dimension of
the final dynamics is related to the verification of
criterion~(GrF). In~[10], a scalar parabolic equation of the form
(1.2) is considered in a rectangle with Dirichlet boundary
condition. The authors present conditions under which the
attractor is a Lipschitz graph over finitely many first modes of
the Laplace operator. And they use the \textit{cone condition}
well-known in the literature~[13,16,19]. In this connection, it
seems to be very perspective to study the problem of
finite-dimensional reduction of systems of equations (1.2) on the
two-dimensional torus~$\mathbb{T}^{2}$.

\medskip
\textbf{Acknowledgements.} I am grateful to Sergey Zelik for
useful discussions.

\bigskip
\centerline{\large\bf References}
\medskip

\noindent [1] A. V. Babin and M. I. Vishik, \textit{Attractors of
evolution equations}, North-Holland, 1992.

\noindent [2] Ju. L. Daleckii and M. G. Krein, \textit{Stability
of solutions of differential equations in}

\textit{Banach space}, Trans. Math. Monogr., vol. \textbf{43},
Amer. Math. Soc., Providence,

Rhode Island, 1974.

\noindent [3] A. Eden, C. Foias, B. Nicolaenko and R. Temam,
\textit{Exponential Attractors for }

 \textit{Dissipative Evolution Equations}, Wiley, New York, 1994.

\noindent [4] A. Yu. Goritskii and V. V. Chepyzhov, ``The
dichotomy property of solutions of

quasilinear equations in problems on inertial manifolds'',
\textit{Sb. Mathematics}, \textbf{196}:3-4

(2005), 485--511.

\noindent [5] D. Henry, \textit{Geometric theory of semilinear
parabolic equations}, Lect. Notes in Math.,

vol. \textbf{840}, Springer, 1981.

\noindent [6] D. A. Kamaev, ``Families of stable manifolds of
invariant sets of systems of parabolic

equations'', \textit{Russ. Math. Surv.}, \textbf{47}:5 (1992),
185--186.

\noindent [7] A. Kostianko and S. Zelik, ``Inertial manifolds for
1D reaction-diffusion-advection

systems. Part I: Dirichlet and Neumann boundary conditions'',
\textit{Comm. Pure Appl.}

\textit{Anal.,} \textbf{16}:6 (2017), 2357--2376.

\noindent [8] A. Kostianko and S. Zelik, ``Inertial manifolds for
1D reaction-diffusion-advection

systems. Part II: Periodic boundary conditions'', \textit{Comm.
Pure Appl. Anal}., \textbf{17}:1

(2018), 285--317.

\noindent [9] I. Kukavica, ``Fourier parametrization of attractors
for dissipative equations in one

space dimension'',  \textit{J. Dyn. Differ. Eq.}, \textbf{15}:2--3
(2003), 473--484.

\noindent [10] M. Kwak and B. Lkhagvasuren, ``The cone property
for a class of parabolic

equations``, \textit{Journal of the Korean SIAM}, \textbf{21}:2
(2017), 81--87.


\noindent [11] P. Lankaster, \textit{Theory of matrices}, Academic
Press, New York -- London, 1969.

\noindent [12] J. C. Robinson, ``Global attractors: topology and
finite-dimensional dynamics'',

\textit{J. Dyn. Differ. Eq.}, \textbf{11}:3 (1999), 557--581.

\noindent [13] J. C. Robinson, \textit{Infinite-Dimensional
Dynamical Systems},  Cambridge Texts in

Applied Mathematics, Cambridge University Press, 2001.

\noindent [14] A. V. Romanov, ``Finite-dimensional limit dynamics
of dissipative parabolic

equations'', \textit{Sb. Mathematics}, \textbf{191}:3 (2000),
415--429.

\noindent  [15] A. V. Romanov, ``Finite-dimensionality of dynamics
on an attractor for non-linear

parabolic equations'', \textit{Izv. Math}., \textbf{65}:5 (2001),
977--1001.

\noindent [16] R. Temam, \textit{Infinite-dimensional Dynamical
Systems in Mechanics and Physics},

Appl. Math. Sci., vol. \textbf{68} (2-nd ed.), Springer, N.Y.,
1997.

\noindent [17] H. Triebel, ``Theory of function spaces'', Monogr.
in Math., Vol. \textbf{78}, Birkhauser

Verlag, Basel--Boston--Stuttgart, 1983.

\noindent [18] V. I. Yudovich, ``Periodic differential equations
with self-adjoint monodromy

operator'', \textit{Sb. Mathematics}, \textbf{192}:3 (2001),
455--478.

\noindent [19] S. Zelik, ``Inertial manifolds and
finite-dimensional reduction for dissipative PDEs,''

\textit{Proc. Roy. Soc. Edinburgh, Ser. A}, \textbf{144}:6 (2014),
1245--1327.

\bigskip
\bigskip
\noindent School of Applied Mathematics

\noindent National Research University Higher School of Economics

\noindent 34 Tallinskaya St., Moscow, 123458 Russia

\noindent E-mail address: av.romanov@hse.ru

\end{document}